\documentclass[11pt]{article}

\usepackage{amsmath}
\usepackage{amssymb}
\usepackage{DGStpp}
\usepackage{pproof}
\usepackage{remexpp}
\usepackage{mathrsfs}
\let\rscr=\mathscr
\let\mathscr=\relax

\usepackage{eucal}
\let\escr=\mathcal
\let\mathcal=\relax

\newtheorem{theorem}{Theorem}[section]
\newtheorem{proposition}[theorem]{Proposition}
\newtheorem{lemma}[theorem]{Lemma}
\newtheorem{corollary}[theorem]{Corollary}
\newremark{definition}[theorem]{Definition}
\newremark{example}[theorem]{Example}
\newremark{remark}[theorem]{Remark}
\newremark{remarks}[theorem]{Remarks}

\newcommand{\bs}{\backslash}
\newcommand{\bul}{\,\scriptscriptstyle\bullet}
\newcommand{\bve}{\bar{\varepsilon}}
\newcommand{\del}[1]{\nabla_{#1}}
\newcommand{\dg}{\dot\gamma}

\newcommand{\ds}{\oplus} 
\newcommand{\dx}{\dot{x}}
\newcommand{\ddx}{\ddot{x}}

\newcommand{\g}{\gamma}

\newcommand{\half}{\mbox{$\txs\frac{1}{2}$}}
\newcommand{\io}{\iota}
\newcommand{\Iaut}{I^{a\kern-.05em u\kern-.02em t}}
\newcommand{\Ispl}{I^{s\kern-.055em p\kern-.025em l}}
\newcommand{\lsp}{[\kern-0.15em[} 

\newcommand{\n}{\mbox{$\mathfrak{n}$}}
\newcommand{\om}{\omega}
\newcommand{\ph}{$p\kern-.1em H\!$}
\newcommand{\quar}{\mbox{$\txs\frac{1}{4}$}}
\newcommand{\rsp}{]\kern-0.15em]} 
\newcommand{\sC}{\escr C}
\newcommand{\sJ}{\escr J}

\newcommand{\surj}{\rightarrow\kern-.82em\rightarrow}
\newcommand{\tphi}{\tilde{\phi}}
\newcommand{\txs}{\textstyle}
\newcommand{\ve}{\mbox{$\varepsilon$}}
\newcommand{\vp}{\mbox{$\varphi$}}

\newcommand{\z}{\mathfrak{z}}
\newcommand{\E}{\mathfrak{E}}

\newcommand{\G}{\mbox{$\Gamma$}}
\newcommand{\I}{\rscr{I}}
\newcommand{\II}{\rscr{I{\kern-.55em}I}}

\newcommand{\R}{\mbox{${\mathbb R}$}}
\newcommand{\U}{\mathfrak{U}}
\newcommand{\V}{\mathfrak{V}}
\newcommand{\Z}{\mathfrak{Z}}

\renewcommand{\v}{\mathfrak v}

\makeatletter
\newcommand{\Ad}[1]{\mathop{\operator@font Ad}\nolimits_{#1}}
\newcommand{\ad}[1]{\mathop{\operator@font ad}\nolimits_{#1}}
\newcommand{\add}[2]{\mathop{\operator@font ad}\nolimits^{\dagger}_{#1}{#2}}
\newcommand{\Aut}{\mathop{\operator@font Aut}\nolimits}
\newcommand{\diag}{\mathop{\operator@font diag}\nolimits}
\newcommand{\End}{\mathop{\operator@font End}\nolimits}
\newcommand{\Ric}{\mathop{\operator@font Ric}\nolimits}
\newcommand{\Dspecp}{\mathop{\operator@font {\escr D}spec}\nolimits_\wp}
\newcommand{\specl}{\mathop{\operator@font spec}\nolimits_\ell}
\newcommand{\specp}{\mathop{\operator@font spec}\nolimits_\wp}
\newcommand{\tr}{\mathop{\operator@font tr}\nolimits}
\makeatother

\font\heads = cmbx12

\preprint{CP6}
\title{Lattices and Periodic Geodesics\\ in Pseudoriemannian 2-step
Nilpotent Lie Groups}
\author{Luis A. Cordero\thanks{Partially supported by Project
        MEC:MTM2005-08757-C04-01, Spain.}}
\address{Dept. Xeometr{\'\i}a e Topolox{\'\i}a\\
        Facultade de Matem\'aticas\\
        Universidade de Santiago\\
        15782 Santiago de Compostela\\
        Spain\\
        cordero@zmat.usc.es}
\author{Phillip E. Parker\thanks{Partially supported by MEC:DGES
Program SAB1995-0757, Spain.}}
\address{Mathematics Department\\
        Wichita State University\\
        Wichita KS 67260-0033\\
        USA\\
        phil@math.wichita.edu}
\date{7 Feb 2008}

\abstract{We give a basic treatment of lattices $\G$ in these groups.
Certain tori $T_F$ and $T_B$ provide the model fiber and the base for a
submersion of $\G\bs N$.  This submersion may not be pseudoriemannian in
the usual sense, because the tori may be degenerate.  We then begin the
study of periodic geodesics in these compact nilmanifolds, obtaining a
complete calculation of the period spectrum of certain flat spaces.

Appeared in {\it Int.\ J.\ Geom.\ Methods Mod.\ Phys.} {\bf 5} (2008)
79--99.
}

\msc{53C50}{22E25, 53B30, 53C30}

\begin{document}

\maketitle


\section{\heads Introduction}

The 2-step nilpotent groups are nonabelian and as close as possible to
being Abelian, but display a rich variety of new and interesting geometric
phenomena.  As in the Riemannian case, one of many places where they arise
naturally is as groups of isometries acting on horospheres in certain
(pseudoriemannian) symmetric spaces.  Another is in the Iwasawa
decomposition of semisimple groups with the Killing metric, which need not
be definite.  Here we study the lattices and periodic geodesics in these
group spaces. For a more extensive historical introduction that better
puts them into the contemporary context, see \cite{CP4,CP5}. A recent,
masterful survey of the Riemannian case is \cite{E2}.

By an {\em inner product\/} on a vector space $V$ we shall mean a
nondegenerate, symmetric bilinear form on $V$, generally denoted by
$\langle\,,\rangle$.  In particular, we {\em do not\/} assume that it is
positive definite.  Our convention is that $v\in V$ is timelike if
$\langle v,v\rangle > 0$, null if $\langle v,v\rangle = 0$, and spacelike
if $\langle v,v\rangle < 0$.

Throughout, $N$ will denote a connected, 2-step nilpotent Lie group with
Lie algebra \n\ having center $\z$.  (Recall that 2-step means
$[\n,\n]\subseteq\z$.) We shall use $\langle\,,\rangle$ to denote either
an inner product on \n\ or the induced left-invariant pseudoriemannian
(indefinite) metric tensor on $N$.

In Section~\ref{defex}, we give the fundamental definitions and examples
used in the rest of this paper.  The main problem encountered is that the
center $\z$ of \n\ may be degenerate:  it might contain a (totally) null
subspace.  We shall see that this possible degeneracy of the center causes
the essential differences between the Riemannian and pseudoriemannian cases.

Much like the Riemannian case, we would expect that $(N,\langle\, ,\rangle)$
should in some sense be similar to flat pseudoeuclidean space.  This is seen
in the examples of totally geodesic subgroups in Section \ref{geods}.  We
also showed \cite{CP4} the existence of $\dim\z$ independent first
integrals, a familiar result in pseudoeuclidean space.  Unlike the
Riemannian case, there are flat groups which are isometric to
pseudoeuclidean spaces.

Section~\ref{latpg} begins with a basic treatment of lattices $\G$ in these
groups.  The tori $T_F$ and $T_B$ provide the model fiber and the base for a
submersion of $\G\bs N$.  This submersion may not be pseudoriemannian in the
usual sense, because the tori may be degenerate.  Also, it is possible for
a (null) geodesic to be closed but not periodic. Thus we begin the study
of periodic geodesics in such a compact nilmanifold and of its period
spectrum.  In the Riemannian case, this is its length spectrum.  We obtain
a complete calculation of the period spectrum for the flat spaces of
Section \ref{defex}.  For related work on the length spectrum in the
Riemannian setting (which is closely related to the geometry of the
Laplacian and which plays a central role in isospectral questions), we
refer to \cite{CGMYDS97, RG00, RG05, L00, MMRG00, MMRG04}.

We recall \cite{CP4,E2} some basic facts about 2-step nilpotent Lie
groups.  As with all nilpotent Lie groups, the exponential map $\exp
: \n\rightarrow N$ is surjective.  Indeed, it is a diffeomorphism for
simply connected $N$; in this case we shall denote the inverse by $\log$.
The Baker-Campbell-Hausdorff formula takes on a particularly simple form
in these groups:
\begin{equation}
\label{bch}
\exp(x)\exp(y)=\exp(x+y+\half[x,y])\, .
\end{equation}
Letting $L_n$ denote left translation by $n\in N$, we have the following:
\begin{lemma}
Let\label{dexp} \n\ denote a 2-step nilpotent Lie algebra and $N$ the
corresponding simply connected Lie group. If $x,a\in\n$, then
$$ \exp_{x*}(a_x)=L_{\exp(x)*}\left( a+\half[a,x]\right) $$
where $a_x$ denotes the initial velocity vector of the curve $t\mapsto
x+ta$.
\end{lemma}
\begin{corollary}
In\label{pcc} a pseudoriemannian 2-step nilpotent Lie group, the
exponential map preserves causal character.  Alternatively, 1-parameter
subgroups are curves of constant causal character.
\end{corollary}
\begin{proof}
For the 1-parameter subgroup $c(t) = \exp(t\,a)$, one easily sees that $\dot{c}(t) =
\exp_{ta*}(a) = L_{\exp(ta)*}a$ and left translations are isometries.
\end{proof}
Of course, 1-parameter subgroups need not be geodesics, as simple
examples show \cite{CP4}.

We shall also need some basic facts about lattices in $N$. In nilpotent
Lie groups, a lattice is a discrete subgroup \G\ such that the homogeneous
space $M=\G\bs N$ is compact \cite{R}. Lattices do not always exist in
nilpotent Lie groups \cite{M}.
\begin{theorem}
The simply connected, nilpotent Lie group $N$ admits a lattice if and only
if there exists a basis of its Lie algebra \n\ for which the structure
constants are rational.
\end{theorem}
Such a group is said to have a rational structure, or simply to be
rational.

Geodesic completeness is notoriously problematic in pseudoriemannian
spaces.  For 2-step nilpotent Lie groups, things work nicely as shown by
this result first published by Guediri \cite{G}.
\begin{theorem}
On\label{gt} a 2-step nilpotent Lie group, all left-invariant pseudoriemannian
metrics are geodesically complete.
\end{theorem}
He also provided an explicit example of an incomplete metric on a 3-step
nilpotent Lie group; we refer to \cite{CP4} for complete calculations of explicit formulas for
the connection, curvatures, covariant derivative, {\em etc.}

\section{\heads Definitions and Examples}
\label{defex}

For the convenience of the reader, we repeat some basic definitions and
theorems from \cite{CP5}. Complete details and proofs are in \cite{CP4}.

In the Riemannian (positive-definite) case, one splits $\n = \z\ds\v =
\z\ds\z^{\perp}$ where the superscript denotes the orthogonal complement
with respect to the inner product $\langle\,,\rangle$.  In the general
pseudoriemannian case, however, $\z\ds\z^{\perp}\not=\n$.  The problem is
that $\z$ might be a degenerate subspace; {\em i.e.,} it might contain a
null subspace $\U$ for which $\U\subseteq\U^{\perp}$.

Thus we have to adopt a more complicated decomposition of \n.
Observe that if $\z$ is degenerate, the null subspace $\U$ is well defined
invariantly. We use a decomposition
$$
\n = \z\ds\v = \U\ds\Z\ds\V\ds\E
$$
in which $\z = \U\ds\Z$ and $\v = \V\ds\E$, $\U$ and $\V$ are complementary
null subspaces, and $\U^{\perp}\cap\V^{\perp} = \Z\ds\E$. Although the
choice of $\V$ is {\em not\/} well defined invariantly, once a $\V$ has been
chosen then $\Z$ and $\E$ {\em are\/} well defined invariantly. Indeed, $\Z$
is the portion of the center $\z$ in $\U^{\perp}\cap\V^{\perp}$ and $\E$
is its orthocomplement in $\U^{\perp}\cap\V^{\perp}$. This is a Witt
decomposition of \n\ given $\U$ as described in \cite[p.\,37f\,]{P}, easily
seen by noting that $\left(\U\ds\V\right)^{\perp} = \Z\ds\E$, adapted to the
special role of the center in \n.

We fix a choice of $\V$ (and therefore $\Z$ and $\E$). Having fixed
$\V$, observe that the inner product $\langle\,,\rangle$ provides
a dual pairing between $\U$ and $\V$; {\em i.e.,} isomorphisms $\U^*\cong\V$
and $\U\cong\V^*$. Thus the choice of a basis $\{u_i\}$ in $\U$ determines an
isomorphism $\U\cong\V$ ({\em via\/} the dual basis $\{v_i\}$ in $\V$).
In addition to the choice of $\V$, we also fix a basis of $\U$.

We also need an involution $\io$ that interchanges $\U$ and
$\V$ by this isomorphism and which reduces to the identity on $\Z\oplus\E$
in the Riemannian (positive-definite) case. The choice of such an
involution is not significant \cite{CP4}.  In terms of chosen orthonormal
bases
$\{z_\alpha\}$ of $\Z$ and $\{e_a\}$ of $\E$,
$$
\io (u_i)=v_i, \quad \io (v_i)=u_i, \quad \io (z_\alpha )=\ve_\alpha \,
z_\alpha , \quad \io (e_a )= \bar{\ve}_a \, e_a\,,
$$
where, as usual,
$$
\langle u_i,v_i \rangle = 1, \quad \langle z_\alpha ,z_\alpha \rangle =
\ve_\alpha, \quad \langle e_a,e_a \rangle = \bar{\ve}_a\,.
$$
Then $\io (\U )=\V$, $\io (\V )=\U$, $\io (\Z )=\Z$, $\io (\E )=\E$ and
$\io^2 = I$.  It is obvious that $\io$ is selfadjoint with respect to the
inner product,
\begin{equation}
\langle \io x,y \rangle = \langle x, \io y \rangle , \qquad x,y \in \n\,,
\label{iosa}
\end{equation}
so $\io$ is an isometry of \n. (However, it does {\em not\/} integrate to
an isometry of $N$; see Example \ref{hq}.) Moreover,
\begin{equation}
\langle x, \io x \rangle = 0 \mbox{ if and only if } x = 0, \qquad x \in \n\,.
\label{iosnd}
\end{equation}

Consider the adjoint with respect to $\langle\,,\rangle$ of the
adjoint representation of the Lie algebra \n\ on itself, denoted by
$\add{}{}\!$.  First note that for all $a\in\z$,
$\add{a}{\!\vcenter{\hbox{$\bul$}}} = 0$.  Thus for all $y\in\n$,
$\add{\bul}{y}$ maps $\V\ds\E$ to $\U\ds\E$.  Moreover, for all $u\in\U$
we have $\add{\bul}{u} = 0$ and for all $e\in\E$ also $\add{\bul}{e} = 0$.
Following \cite{K,E,E1}, we define the operator $j$. Note the use of
the involution $\io$ to obtain a good analogy to the Riemannian case.
\begin{definition}
The\label{dj} linear mapping
$j:\U\ds\Z\rightarrow\End\left(\V\ds\E\right) $
is given by
$j(a)x = \io\add{x}{\io a}$.
\end{definition}

Let $x,y\in\n$. Recall \cite{BP2,CP4} that homaloidal planes are those for
which the numerator $\langle R(x,y)y,x\rangle$ of the sectional curvature
$K(x,y)$ vanishes.  This notion is useful for degenerate planes tangent to
spaces that are not of constant curvature.
\begin{theorem}
All\label{zph} central planes are homaloidal: $R(z,z')z'' = R(u,x)y = R(x,y)u
= 0$ for all $z,z',z''\in\Z$, $u\in\U$, and $x,y\in\n$. Thus the
nondegenerate part of the center is flat:
$K(z,z')=0$.
\end{theorem}
In view of this result, we extend the notion of flatness to possibly
degenerate submanifolds.
\begin{definition}
A submanifold of a pseudoriemannian manifold is {\em flat\/} if and only
if every plane tangent to the submanifold is homaloidal.
\end{definition}
\begin{corollary}\label{zf}
The center $Z$ of $N$ is flat.
\end{corollary}
\begin{corollary}\label{ccf}
The only $N$ of constant curvature are flat.
\end{corollary}
The degenerate part of the center can have a profound effect on the
geometry of the whole group.
\begin{theorem}
If\/ $[\n,\n] \subseteq \U$ and\label{e0f} $\E=\{0\}$, then $N$ is flat.
\end{theorem}

We continue with a formula for certain sectional curvatures.
\begin{theorem}
If\label{kee} $e,e'$ are any orthonormal vectors in $\E$, then
$$
K(e,e')=-{\txs\frac{3}{4}}\bve\bve{\kern.05em}'\langle [e,e'],[e,e']\rangle
$$
with $\bve=\langle e,e\rangle$ and $\bve{\kern.05em}'=\langle e',e'
\rangle$.
\end{theorem}
Some sectional curvature numerators are also relevant.
\begin{proposition}
If\label{oscn} $z\in\Z$, $v\in\V$, and $e\in\E$, then
\begin{eqnarray*}
\langle R(z,v)v,z\rangle &=& \quar\langle j(\io z)v,j(\io z)v\rangle ,\\
\langle R(v,e)e,v\rangle &=& -{\txs\frac{3}{4}}\langle [v,e],[v,e]\rangle
 +\quar\langle j(\io v)e,j(\io v)e\rangle ,\\
\langle R(v,v')v',v\rangle &=& -{\txs\frac{3}{4}}\langle [v,v'],[v,v']\rangle
 +\half\langle j(\io v)v',j(\io v')v\rangle \\
 & & {}+\quar\Bigl( \langle j(\io v')v,j(\io v')v\rangle
     +\langle j(\io v)v',j(\io v)v'\rangle\Bigr)\\
 & & {}-\langle j(\io v)v,j(\io v')v'\rangle .
\end{eqnarray*}
\end{proposition}

Here is the example we mentioned previously; details and other examples
are in \cite{CP4,CP5}.
\begin{example}
For\label{hq} the simplest quaternionic Heisenberg algebra of dimension 7,
we may take a basis $\{u_1,u_2,z,v_1,v_2,e_1,e_2\}$ with structure
equations
$$ \begin{array}{rclcrcl}
[e_1,e_2]&=&z &\qquad & [v_1,v_2]&=&z\\
{[}e_1,v_1]&=&u_1 && [e_2,v_1]&=&u_2\\
{[}e_1,v_2]&=&u_2 && [e_2,v_2]&=&-u_1
\end{array} $$
and nontrivial inner products
$$ \langle u_i,v_j\rangle=\delta_{ij}\, ,\quad\langle z,z\rangle=\ve\,
,\quad\langle e_a,e_a\rangle=\bve_a\,. $$
As usual, each \ve-symbol is $\pm 1$ independently (this is a combined
null and orthonormal basis), so the signature is
$(++--\,\ve\,\bve_1\,\bve_2)$.

Sectional curvatures for this group are
\begin{eqnarray*}
\langle R(v_1,v_2)v_2,v_1\rangle&=&-(\bve_1+{\txs\frac{3}{4}}\ve)\,,\\
\langle R(v,e)e,v\rangle=\langle R(z,v)v,z\rangle &=& 0\, ,\\
K(z,e_1)=K(z,e_2) &=& \quar\ve\bve_1\bve_2\, ,\\
K(e_1,e_2) &=& -{\txs\frac{3}{4}}\ve\bve_1\bve_2\, .
\end{eqnarray*}
Thus $\io$ cannot integrate to an isometry of $N$ in general, as mentioned
after equation (\ref{iosa}). Isometries must preserve vanishing of
sectional curvature, and an integral of $\io$ would interchange homaloidal
and nonhomaloidal planes in this example.
\end{example}

\section{\heads Totally Geodesic Subgroups and Geodesics}
\label{geods}

We begin by noting that O'Neill \cite[Ex.\,9, p.\,125]{O} has extended the
definition of totally geodesic to degenerate submanifolds of
pseudoriemannian manifolds. We shall use this extended version.

Recall from \cite{O} that the extrinsic and intrinsic curvatures of
totally geo\-des\-ic submanifolds coincide. Thus there is an unambiguous
notion of flatness for them.  Note that a connected subgroup $N'$ of $N$
is a totally geodesic submanifold if and only if it is totally geodesic at
the identity element of $N$, because left translations by elements of $N'$
are isometries of $N$ that leave $N'$ invariant. A connected, totally
geodesic submanifold need not be a connected, totally geodesic subgroup,
but $(N,\langle\, ,\rangle)$ has many totally geodesic subgroups. Many of
them are flat, illustrating the similarity to pseudoeuclidean spaces; {\em
cf.} \cite[(2.11)]{E1}.
\begin{example}
For\label{1par} any $x\in\n$ the 1-parameter subgroup $\exp (tx)$ is a
geo\-des\-ic if and only if $\add{x}{x}=0$. We find this if and only if
$x\in\z$ or $x\in\U\ds\E$. This is essentially the same as the Riemannian
case, but with some additional geodesic 1-parameter subgroups coming from
$\U$.
\end{example}
\begin{example}
Abelian\label{abel} subspaces of $\V\ds\E$ are Lie subalgebras of \n, and
give rise to complete, flat, totally geodesic abelian subgroups of $N$,
just as in the Riemannian case \cite{E1}. The construction given in
\cite[(2.11),\,Ex.\,2]{E1} is valid in general, and shows that if $\dim\V
\ds\E \geq 1 + k + k\dim\z$, then every nonzero element of $\V\ds\E$ lies
in an abelian subspace of dimension $k+1$.
\end{example}
\begin{example}
The\label{ftgZ} center $Z$ of $N$ is a complete, flat, totally geodesic
submanifold.  Moreover, it determines a foliation of $N$ by its left
translates, so each leaf is flat and totally geodesic, as in the
Riemannian case \cite{E1}. In the pseudoriemannian case, this foliation in
turn is the orthogonal direct sum of two foliations determined by $\U$ and
$\Z$, and the leaves of the $\U$-foliation are also null.  All these
leaves are complete.
\end{example}

For the geodesic equation, let us consider (as suffices) a geodesic $\g$
with $\g(0) = 1\in N$ and $\dot{\g}(0) = a_0+x_0\in \z\ds\v$. One may further
decompose $a_0 = u_0 + z_0 \in \U\ds\Z$ and $x_0 = v_0 + e_0 \in \V\ds\E$.
In exponential coordinates, write $\g(t) = \exp \left( a(t)+x(t)\right) =
\exp\left( u(t) + z(t) + v(t) + e(t)\right)$ where one has $\dot{a}(0) = a_0$,
$\dot{x}(0) = x_0$, {\em etc.} For the tangent vector $\dot{\g}$, we obtain
\begin{eqnarray*}
\dg &=& \exp_{(a+x)*}(\dot{a}+\dx)\\
   &=& L_{\g(t)*}\left(\dot{a}+\dx + \half[\dot{a}+\dx,a+x]\right)\\
   &=& L_{\g(t)*}\left(\dot{a}+\dx + \half[\dx,x]\right),
\end{eqnarray*}
using Lemma \ref{dexp},
regarded as vector fields along $\g$. Then the geodesic equation is
equivalent to
$$ \frac{d}{dt}\left( \dot{a}+\dx + \half[\dx,x]\right)
   - \add{\dot{a} + \dx +\frac{1}{2}[\dx,x]}{\left(\dot{a}
   + \dx +\half[\dx,x]\right)} = 0\,. $$
Simplifying slightly, we find
\begin{equation}
\underbrace{\frac{d}{dt}\left(\dot{a} + \half[\dx,x]
\right)_{\strut}}_{\in\,\U\ds\Z}{\;} +
\underbrace{\ddx_{\strut}}_{\in\,\V\ds\E} -
{\;\,}\underbrace{\add{\dx}{\left(\dot{z} + \dot{v} +
\half[\dx,x]\right)}_{\strut}}_{\in\,\U\ds\E} = 0\,.\label{ges}
\end{equation}
Using superscripts to denote components, we obtain for the
$\Z$-component
$$ \frac{d}{dt}\left(\dot{z} + \half[\dx,x]^{\Z}\right) = 0\,. $$
Using the initial condition, we get
$$ \dot{z} + \half[\dx,x]^{\Z} = z_0\,. $$
Next we note that $\ddot{v} = 0$ whence $\dot{v}(t) = v_0$ is a
constant. We use these to simplify the other component equations.
\begin{eqnarray}
\frac{d}{dt}\left(\dot{u} + \half[\dx,x]^{\U}\right) -
   \left(\add{\dx}{(z_0+v_0)}\right)^{\U} &=& 0\label{ges1}\\
\dot{z} + \half[\dx,x]^{\Z} &=& z_0\label{ges2}\\
\ddot{v} &=& 0\label{ges3}\\
\ddot{e} - \left(\add{\dx}{(z_0+v_0)}\right)^{\E} &=& 0\label{ges4}
\end{eqnarray}

In analogy with Eberlein \cite{E,E1} we define two operators.
\begin{definition}
For\label{defJ} fixed $z_0\in\Z$ and $v_0\in\V$ as above, define
\begin{eqnarray*}
\rscr{J}:\V\ds\E \longrightarrow \E &:& y \longmapsto
\left(\add{y}{(z_0+v_0)}\right)^{\E},\\
\sJ:\V\ds\E \longrightarrow \U &:& y \longmapsto
\left(\add{y}{(z_0+v_0)}\right)^{\U}.
\end{eqnarray*}
\end{definition}
We shall denote the restriction of $\rscr{J}$ to $\E$ by $J$, and this will
play the same role as $J$ in Eberlein \cite{E,E1}.

Now we rewrite the geodesic equations in terms of $J$ and $\sJ$,
using the linearity of $\add{}{}$ to rearrange some terms.
\begin{eqnarray}
\frac{d}{dt}(\dot{u} + \half[\dx,x]^{\U}) &=& \sJ\dx \label{ge1}\\
\dot{z} + \half[\dx,x]^{\Z} &=& z_0\label{ge2}\\
\ddot{v} &=& 0\label{ge3}\\
\ddot{e} - J\dot{e} &=& \rscr{J}v_0\label{ge4}
\end{eqnarray}
While the $\V$-component of a geodesic is simple, its mere presence
affects all of the other components.

We also readily see that the system is completely integrable.  Thus, as
noted in the Introduction (Theorem \ref{gt}), all left-invariant
pseudoriemannian metrics on these groups are complete.  Also, {\em
regardless of signature,} we may obtain the existence of $\dim\z$ first
integrals as in \cite{E1}.

Keep $J\in\End(\E)$ and observe that $J$ is skewadjoint with respect to
$\langle\,,\rangle$.  Write $\E = \E_1\ds\E_2$ with $\E_1 = \ker J$ as a
direct sum.  Unfortunately, it need not be orthogonal; see \cite{JP2} for
a complete list of canonical forms for such $J$.  Thus we shall assume for
the rest of this article that $\E_1\ds\E_2$ is orthogonal.  For example,
this will be the case if $\ker J$ is nondegenerate.

Decompose $\rscr{J}v_0 = y_1 + y_2 \in \E_1\ds\E_2$, respectively.  Note
that $J$ is invertible on $\E_2$; we denote this restriction by $J$ also.

We continue with a geodesic through the identity element, so $\g(0)=1$ and
$\dot{\g}(0) = a_0+x_0 = u_0+z_0+v_0+e_0$, with $J\in\End(\E)$ and $\E =
\E_1\ds\E_2$ with $\E_1 = \ker J$, and with $\rscr{J}v_0 = y_1 + y_2 \in
\E_1\ds\E_2$ as before.  Decompose $e_0 = e_1 + e_2 \in \E_1\ds\E_2$,
respectively.  (These $e_i$ should not be confused with the basis elements
appearing in other sections.)  Recall that $J$ is invertible on $\E_2$,
and that we let $J$ denote the restriction there as well.  For
convenience, set
 $x_1= e_1 + v_0 - J^{-1}y_2$ and
$x_2= e_2 + J^{-1}y_2$.
\begin{theorem}
Under\label{ige} these assumptions, the geodesic equations may be
integrated as:
\begin{eqnarray}
x(t) &=& t\,x_1 + \left( e^{tJ} - I\right) J^{-1}x_2 + \half t^2y_1\,,
         \label{ige1}\\
z(t) &=& t\,z_0 + \I[\dx,x]^{\Z},\label{ige2}\\
u(t) &=& t\,u_0 + \I[\dx,x]^{\U} + \II\sJ\dx\,,\label{ige3}
\end{eqnarray}
where
\begin{eqnarray*}
 \I[\dx,x] &=& -\half\int_0^t\left[\dx(s),x(s)\right] ds\\
&=&\half t\left[ x_1 + \half ty_1, \left(e^{tJ} + I\right) J^{-1}x_2 \right]
- t\left[ y_1, e^{tJ}x_2\right] + {\txs\frac{1}{12}}t^3[x_1, y_1]\\
&&{}+ \half\left[\left( e^{tJ} - I\right) J^{-1}x_2, J^{-1}x_2\right]
- \left[ x_1, \left( e^{tJ} - I\right) J^{-2}x_2\right]\\
&&{}+ \left[ y_1, \left( e^{tJ} - I\right) J^{-3}x_2\right]
+ \half\int_0^t \left[ e^{sJ} J^{-1}x_2, e^{sJ} x_2\right] ds
\end{eqnarray*}
and
$$ \II\sJ\dx = \int_0^t\int_0^s\sJ\dx(\sigma)\,d\sigma\,ds\,. $$
\end{theorem}
\begin{proof}
The formulas follow from straightforward integrations of the geodesic
equations (\ref{ge1})--(\ref{ge4}).  We used the general fact about
exponentials of matrices that $J$ commutes with $e^{tJ}$ for all $t\in\R$.
Using this, it is routine to verify that $x(t)$, $z(t)$, and $u(t)$
satisfy the geodesic equations and initial conditions.
\end{proof}
\begin{corollary}
When\label{igend} $N$ has a nondegenerate center, the formulas simplify
somewhat.  Now equation (\ref{ge4}) is homogeneous and we obtain
\begin{eqnarray}
e(t) &=& t\,e_1 + \left( e^{tJ} - I\right) J^{-1}e_2\,,\label{igend1}\\
z(t) &=& t\,z_1(t) + z_2(t) + z_3(t)\,,\label{igend2}
\end{eqnarray}
where
\begin{eqnarray*}
z_1(t) &=& z_0+\half\left[e_1,\left(e^{tJ}+I\right)J^{-1}e_2\right],\\
z_2(t) &=& \left[e_1,\left(I-e^{tJ}\right) J^{-2}e_2\right] +
           \half\left[e^{tJ}J^{-1}e_2, J^{-1}e_2\right],\\
z_3(t) &=& \half\int_0^t \left[e^{sJ}J^{-1}e_2, e^{sJ}e_2\right]
      \,ds\,.
\end{eqnarray*}
Note that $z_3$ may contribute to $z_1$ and $z_2$.\eop
\end{corollary}
The flat spaces we found in Theorem \ref{e0f} allow more simplification.
\begin{corollary}
When\label{iged} $[\n,\n] \subseteq \U$ and $\E = \{0\}$, then $x = v$ and
we obtain
\begin{eqnarray}
v(t) &=& t\,v_0\,, \label{iged1}\\
z(t) &=& t\,z_0\,, \label{iged2}\\
u(t) &=& t\,u_0 + \half t^2 \,\sJ v_0\,. \label{iged3}
\end{eqnarray}
\end{corollary}
\begin{corollary}
If\label{dgc} $[\n,\n] \subseteq \U$ and $\E=\{0\}$, then $N$ is
geodesically connected.  Consequently, so is any nilmanifold with such a
universal covering space.
\end{corollary}
\begin{proof}
$N$ is complete by Theorem \ref{gt}, hence pseudoconvex.  The preceding
geodesic equations show that $N$ is nonreturning.  Thus the space of
geodesics $G(N)$ is Hausdorff by Theorem 5.2 of \cite{BP9}.  Now Theorem
4.2 of \cite{BLP} yields geodesic connectedness of $N$.
\end{proof}
Thus these compact nilmanifolds are much like tori.  This is also
illustrated by the computation of their period spectrum in Theorem
\ref{upd}.

\section{\heads Lattices and Periodic Geodesics}
\label{latpg}

In this section, we assume that $N$ is rational and let \G\ be a lattice
in $N$. Since $N$ is 2-step nilpotent, it has nice generating sets of \G;
{\em cf.}~\cite[(5.3)]{E}.
\begin{proposition}
If $N$ is a simply connected, 2-step nilpotent Lie group of dimension $n$
with lattice \G\ and with center $Z$ of dimension $m$, then there exists a\/
{\em canonical} generating set $\{\vp_1,\ldots,\vp_n\}$ such that
$\{\vp_1,\ldots,\vp_m\}$ generate $\G\cap Z$. In particular, $\G\cap Z$
is a lattice in $Z$.
\end{proposition}
From this, formula (\ref{bch}), and \cite[Prop.\,2.17]{R}, one obtains
as in \cite[(5.3)]{E1}
\begin{corollary}
Let\label{vl} $N$ be a simply connected, 2-step nilpotent Lie group with
lattice \G\ and let $\pi:\n\rightarrow\v$ denote the projection. Then
 \begin{enumerate}
  \item $\log\G\cap\z$ is a vector lattice in $\z$;
  \item $\pi(\log\G)$ is a vector lattice in $\v$;
  \item $\G\cap Z=Z(\G)$.\eop
 \end{enumerate}
\end{corollary}
Here, we used the splitting $\n =\U\ds\Z\ds\V\ds\E$ from Section~\ref{defex}
with $\z = \U\ds\Z$ and $\v = \V\ds\E$.
Thus to the compact nilmanifold $\G\bs N$ we may associate two flat
(possibly degenerate) tori; {\em cf.} \cite[p.\,644]{E1}.
\begin{definition}
With\label{tfb} notation as preceding,
\begin{eqnarray*}
T_{\z} &=& \z/(\log\G\cap\z)\, ,\\
T_{\v} &=& \v/\pi(\log\G)\, .
\end{eqnarray*}
\end{definition}
Observe that $\dim T_{\z} + \dim T_{\v} = \dim\z + \dim\v = \dim\n$.

Next, we apply Theorem 3 from \cite{PS} to our situation. Let $T^m$ denote
the $m$-torus as usual.
\begin{theorem}
Let\label{prin} $m=\dim\z$ and $n=\dim\v$. Then $\G\bs N$ is a principal
$T^m$-bundle over $T^n$.
\end{theorem}
The model fiber $T^m$ can be given a geometric structure from its closed
embedding in $\G\bs N$; we denote this geometric $m$-torus by $T_F$.
Similarly, we wish to provide the base $n$-torus with a geometric structure
so that the projection $p_B:\G\bs N\surj T_B$ is the appropriate
generalization of a pseudoriemannian submersion \cite{O} to (possibly)
degenerate spaces. Observe that the splitting $\n = \z\ds\v$ induces
splittings $TN = \z N\ds\v N$ and $T(\G\bs N) =
\z(\G\bs N)\ds\v(\G\bs N)$, and that $p_{B*}$ just mods out $\z(\G\bs N)$.
Examining the definition on page 212 of \cite{O}, we see
that the key is to construct the geometry of $T_B$ by defining
\begin{equation}\label{dprs}
p_{B*} : \v_\eta(\G\bs N)\to T_{p_B(\eta)}(T_B) \mbox{ for each }
\eta\in\G\bs N \mbox{ is an isometry }
\end{equation}
and
\begin{equation}\label{dprc}
\nabla^{T_B}_{p_{B*}x}p_{B*}y = p_{B*}\left(\pi\nabla_x y\right)
\mbox{ for all $x,y\in\v=\V\ds\E$, }
\end{equation}
where $\pi:\n\to\v$ is the projection.  Then the rest of the results of
pages 212--213 in \cite{O} will continue to hold, provided that sectional
curvature is replaced by the numerator of the sectional curvature formula
in his 47.\,Theorem at least when elements of $\V$ are involved:
\begin{equation}\label{o47}
\langle R_{T_B}(p_{B*}x,p_{B*}y)p_{B*}y,p_{B*}x\rangle = \langle
R_{\Gamma\bs N}(x,y)y,x\rangle + {\txs\frac{3}{4}}
\langle[x,y],[x,y]\rangle .
\end{equation}
Now $p_B$ will be a pseudoriemannian submersion in the usual sense if and
only if $\U=\V=\{0\}$, as is always the case for Riemannian spaces.

In the Riemannian case, Eberlein showed that $T_F\cong T_{\z}$ and
$T_B\cong T_{\v}$. Observe that it follows from Theorem \ref{kee},
Proposition \ref{oscn}, and (\ref{o47}) that $T_B$ is flat in general
only if $N$ has a nondegenerate center or is flat.
\begin{proposition}
If\label{scTB} $v,v' \in\V$ and $e,e' \in\E$, then
\begin{eqnarray*}
\langle R_{T_B}(v,e)e,v\rangle &=& \quar\langle j(\io v)e,
        j(\io v)e\rangle  ,\\
\langle R_{T_B}(v,v')v',v\rangle &=& \half\langle j(\io v)v',j(\io
        v')v\rangle - \langle j(\io v)v,j(\io v')v'\rangle\\
 & & {}+\quar\Bigl( \langle j(\io v')v,j(\io v')v\rangle
     +\langle j(\io v)v',j(\io v)v'\rangle\Bigr),\\
K_{T_B}(e,e') &=& 0\,.
\end{eqnarray*}
Here we have suppressed $p_{B*}$ on the left-hand sides for simplicity.
\end{proposition}
Our quaternionic Heisenberg group (Example \ref{hq}) provides a simple
example of a $T_B$ that is not flat.  Indeed, all sectional curvature
numerators vanish except $\langle R_{T_B}(v_1,v_2)v_2,v_1\rangle =
-\bve_1$, where we have again suppressed $p_{B*}$ for simplicity.
\begin{remark}
Observe\label{sTB} that the torus $T_B$ may be decomposed into a topological
product $T_E \times T_V$ in the obvious way.  It is easy to check that
$T_E$ is flat and isometric to $(\log\G \cap \E)\bs\E$, and that $T_V$ has
a linear connection not coming from a metric and not flat in general.
Moreover, the geometry of the product is ``twisted" in a certain way.  It
would be interesting to determine which tori could appear as such a $T_V$
and how.
\end{remark}

We wish to show that the geometry of the fibers $T_F$ is that of $T_{\z}$.
Thus we now consider the submersion $\G\bs N\surj T_B$ and realize the
model fiber $T^m$ as $(\G\cap Z)\bs Z$ considered as tori only, without
geometry.
\begin{definition}
Let\label{defF} $p_N:N\to\G\bs N$ and $p_Z:Z\to T^m = (\G\cap Z)\bs Z$ be the
natural projections.
Define $F:T^m\to I(\G\bs N)$ by $$F\left(p_Z(z)\right)\left(p_N(n)\right) =
p_N(zn)\quad\forall\quad z\in Z,n\in N\,.$$
\end{definition}
We recall from \cite{CP5} that $\Ispl(N)$ denotes the subgroup of the
isometry group $I(N)$ which preserves the splitting $TN = \z N\ds \v N$.
\begin{proposition}
$F$ is\label{Fiso} a smooth isomorphism of groups $T^m\cong\Ispl_0(\G\bs
N)$, where the subscript 0 denotes the identity component.
\end{proposition}
\begin{proof}
We follow Eberlein \cite[(5.4)]{E1}.  It is easy to check that $F$ is a
well-defined, smooth, injective homomorphism with image in $\Ispl_0(\G\bs
N)$. Thus we need only show that $F$ is surjective.

Let $\psi\in\Ispl_0(\G\bs N)$ and let $\phi_t$ be a path from $1 = \phi_0$
to $\psi=\phi_1$. The covering map $p_N$ has the homotopy lifting
property, so choose a lifting $\tphi_t$ as a path in $\Ispl_0(N)$. Then
for all $g\in\G$, it follows that $p_N\left(\tphi_t L_g
\tphi_t^{-1}\right) = p_N$ for all $t$. Hence for each $t\in[0,1]$,
there exists $g_t\in\G$ such that $\tphi_t L_g\tphi_t^{-1} = L_{g_t}$.
Since $g_0 = g$ and $\G$ is a discrete group, it follows that $g_t=g$ for
every $t$, so $L_g$ commutes with $\tphi_t$ for every $g\in\G$ and
$t\in[0,1]$.

From Proposition 3.3 in \cite{CP5}, there exist $n_t\in N$ and $a_t\in
O(N)$ such that $\tphi_t = L_{n_t}a_t$ for all $t$.  Now, every $L_g$
commutes with every $\tphi_t$, such that $a_t(g) = n_t^{-1}gn_t$ for all $t$ and
$g$.  Extension from lattices is unique \cite[Thm.\,2.11]{R}, so $a_t =
\Ad{n_t^{-1}}$.  By Lemma 3.1 in \cite{CP5}, $a_t$ is the identity and
$n_t\in Z$ for all $t$.  Thus $\tphi_1 = L_{n_1}$, so from the definition
of $\tphi_t$ we obtain $p_N L_{n_1} = p_N\tphi_1 = \phi_1 p_N = \psi p_N$.
But this means $F\left(p_Z(n_1)\right) = \psi$.
\end{proof}
\begin{corollary}
$\Ispl_0(\G\bs N)$ acts\label{fftgo} freely on $\G\bs N$ with complete, flat,
totally geodesic orbits.
\end{corollary}
\begin{proof}
By Theorem \ref{prin}, we may identify $\Ispl_0(\G\bs N)$ as the group of
the principal bundle $\G\bs N\surj T_B$, so it acts freely on the total
space.

Since $p_N$ is a local isometry and the $Z$-orbits in $N$ are complete,
flat, and totally geodesic from Example \ref{ftgZ}, it follows (using the
identification $T^m = (\G\cap Z)\bs Z$ {\em supra}) that the
$\Ispl_0(\G\bs N)$-orbits are complete, flat, and totally geodesic.
\end{proof}
\begin{theorem}
Let\label{TF} $N$ be a simply connected, 2-step nilpotent Lie group with
lattice \G, a left-invariant metric tensor, and tori as in the discussion
following Theorem \ref{prin}.  The fibers $T_F$ of the (generalized)
pseudoriemannian submersion $\G\bs N\surj T_B$ are isometric to $T_{\z}$.
\end{theorem}
\begin{proof}
We follow the proof of Eberlein \cite[(5.5),\,item\,2]{E1}. For each $n\in
N$, define $\psi_n = p_N L_n\exp : \z\to\G\bs N$, and note that it is a
local isometry. Clearly, $\psi_n(z) = \psi_n(z')$ if and only if $z' =
z+\log g$ for some $g\in\G\cap Z$. Hence $\psi_n$ induces an isometric
embedding $\tilde{\psi}_n:T_{\z}\to\G\bs N$. That the image is the
$\Ispl_0(\G\bs N)$-orbit of $p_N(n)$ follows from the proof of Corollary
\ref{fftgo}.
\end{proof}
\begin{corollary}
If\label{fTB} in addition the center $Z$ of $N$ is nondegenerate, then
$T_{B}$ is isometric to $T_{\v}$.\eop
\end{corollary}
The proof is essentially the same as the appropriate parts of the proof of
(5.5) in \cite{E1} and we omit it.

We recall that elements of $N$ can be identified with elements of the
isometry group $I(N)$: namely, $n\in N$ is identified with the isometry
$\phi = L_n$ of left translation by $n$. We shall abbreviate this by
writing $\phi\in N$.
\begin{definition}
We\label{dtrl} say that $\phi\in N$ {\em translates\/} the geodesic $\g$ 
by
$\om$ if and only if $\phi\g(t) = \g(t+\om)$ for all $t$. If $\g$ is a
unit-speed geodesic, we say that $\om$ is a {\em period\/} of $\phi$.
\end{definition}
Recall that unit speed means that $|\dot{\g}| = \left|\langle\dot{\g},
\dot{\g}\rangle\right|^{\frac{1}{2}} = 1$.
Since there is no natural normalization for null geodesics, we do not
define periods for them. In the Riemannian case and in the timelike
Lorentzian case in strongly causal spacetimes \cite{BE}, unit-speed
geodesics are parameterized by arclength and this period is a translation
distance. If $\phi$ belongs to a lattice $\G$, it is the length of a
closed geodesic in $\G\bs N$.
\begin{remark}
Note\label{rcc} that it follows from Corollary \ref{pcc} that if $\phi =
\exp(a^* + x^*)$ translates a geodesic $\g$ with $\g(0) = 1 \in N$, then
$a^* + x^*$ and $\dg(0)$ are of the same causal character.
\end{remark}

In general, recall that if $\g$ is a geodesic in $N$ and if $p_N:N\surj
\G\bs N$ denotes the natural projection, then $p_N\g$ is a periodic
geodesic in $\G\bs N$ if and only if some $\phi\in\G$ translates $\g$.
We say {\em periodic\/} rather than {\em closed\/} here because in
pseudoriemannian spaces it is possible for a null geodesic to be closed
but not periodic. If the space is geodesically complete or Riemannian,
however, then this does not occur ({\em cf.}~\cite{O}, p.\,193); the
former is in fact the case for our 2-step nilpotent Lie groups by
Theorem \ref{gt}.  Further recall that free homotopy classes of closed
curves in $\G\bs N$ correspond bijectively with conjugacy classes in $\G$.
\begin{definition}
Let\label{pC} $\sC$ denote either a nontrivial, free homotopy class of
closed curves in $\G\bs N$ or the corresponding conjugacy class in
$\G$. We define $\wp(\sC)$ to be the set of all periods of periodic
unit-speed geodesics that belong to $\sC$.
\end{definition}
In the Riemannian case, this is the set of lengths of closed geodesics in
$\sC$, frequently denoted by $\ell(\sC)$.
\begin{definition}
The\label{psp} {\em period spectrum\/} of $\G\bs N$ is the set
$$ \specp(\G\bs N) = \bigcup_{\sC}\wp(\sC)\,,$$
where the union is taken over all nontrivial, free homotopy
classes of closed curves in $\G\bs N$.
Note that
In the Riemannian case, this is the length spectrum $\specl(\G\bs N)$.
\end{definition}
\begin{example}
Similar\label{sft} to the Riemannian case, we can compute the period
spectrum of a flat torus $\G\bs\R^m$, where $\G$ is a lattice (of maximal
rank, isomorphic to ${\mathbb Z}^m$).  Using calculations related to those of
\cite[pp.\,146--8]{BGM} in an analogous way as for finding the length
spectrum of a Riemannian flat torus, we easily obtain
$$ \specp(\G\bs\R^m) = \{ |g|\ne 0 \mid g\in\G\}\,. $$
It is also easy to see that the nonzero d'Alembertian spectrum is related
to the analogous set produced from the dual lattice $\G^*$ as multiples by
$\pm 4\pi^2$, almost as in the Riemannian case.
\end{example}

As in this example, simple determinacy of periods of unit-speed geodesics
helps make calculation of the period spectrum possible purely in terms of
$\log\G \subseteq \n$.  (See Theorem \ref{upd} for another example.)  Thus
we begin with the following observation.
\begin{proposition}
Let\label{vup} $\phi = \exp(a^* + v^* + e^*)$ translate the unit-speed
geodesic $\g$ by $\om > 0$.  If $v^* \ne 0$, then the period $\om$ is
simply determined.
\end{proposition}
\begin{proof}
We may assume $\g(0) = 1 \in N$ and $\dg(0) = a_0 + v_0 + e_0$.
From Theorem \ref{ige}, $v^* = v(\om) = \om\,v_0$.
\end{proof}
In attempting to calculate the period spectrum then, we can focus our
attention on those cases where the $\V$-component is zero.

From now on, we assume that $N$ is a simply connected, 2-step nilpotent
Lie group with left-invariant pseudoriemannian metric tensor
$\langle\,,\rangle$.  Note that non-null geodesics may be taken to be of
unit speed.  Most nonidentity elements of $N$ translate some geodesic, but
not necessarily one of unit speed; {\em cf.} \cite[(4.2)]{E1}.
\begin{proposition}
Let\label{tsg} $N$ be a simply connected, 2-step nilpotent Lie group with
left-invariant metric tensor $\langle\,,\rangle$ and $\phi\in N$ not the
identity. Write $\log\phi = a^* +x^*\in\z\ds\v$ and assume that $x^* \perp
[x^*,\n]$.  Let $a'$ be the component of $a^*$ orthogonal to $[x^*,\n]$ in
$\z$ and choose $\xi\in\n$ such that $a'=a^* + [x^*,\xi]$.  Set $\om^*
= |a' +x^*|$ if $a'+x^*$ is not null and set $\om^* = 1$ otherwise.  Then
$\phi$ translates the geodesic
$$ \g(t)= \exp(\xi)\,\exp\left(\frac{t}{\om^*}(a'+x^*)\right) $$
by $\om^*$, and $\g$ is of unit speed if $a'+x^*$ is not null.
\end{proposition}
\begin{proof}
Let $n=\exp(\xi)$. One may set $\phi^*=n^{-1}\phi n=\exp(a'+x^*)$ and one may set $\g^*(t)
= n^{-1}\g(t)$. Then $\phi\g(t) = \g(t+\om^*)$ is equivalent to
$\phi^*\g^*(t) = \g^*(t+\om^*)$, and the latter is routine to verify using
(\ref{bch}).
Now $\g$ is a geodesic if and only if $\g^*$ is, and it is easy to check
directly that $\del{a'+x^*}(a'+x^*) = 0$ is equivalent to $\langle
a'+x^*, [x^*,\n]\rangle = 0$.
\end{proof}
Note that the $\U$ components of $a^*$ and $a'$ in fact coincide.  Also
note that if we further decompose $x^* = v^* + e^*$, the result applies to
every nonidentity element $\phi$ with $v^* = 0$.  In particular, when the
center is nondegenerate this is every nonidentity element.
\begin{corollary}
When $\n$ is nonsingular and $\phi\notin Z$, we may take $a'=0$ in
Proposition \ref{tsg}.
\end{corollary}
\begin{proof}
Because then $a^*\in[x^*,\n]$ and $x^*\ne 0$.
\end{proof}

Now we give some general criteria for an element $\phi$ to translate a
geodesic $\g$; {\em cf.} \cite[(4.3)]{E1}.  We use a $J$ as in the passage
following Definition \ref{defJ}, and $x_1$, $x_2$, $y_1$, and $y_2$ as
given just before Theorem \ref{ige}.
\begin{proposition}
Let\label{gctg} $\phi\in N$ and write $\phi = \exp(a^* + x^*)$ for suitable
elements $a^*\in\z$ and $x^*\in\v$. Let $\g$ be a geodesic with $\g(0) =
n\in N$ and $\g(\om) = \phi n$, and let $\dg(0) = L_{n*}(a_0 + x_0)$ for
suitable elements $a_0\in\z$ and $x_0 = v_0+e_0\in\v$.  Let $n^{-1}\g(t)=
\exp\left(a(t) + x(t)\right)$ where $a(t)\in\z$ and $x(t)\in\v$ for all
$t\in\R$ and $a(0) = x(0) = 0$. Then the following are equivalent:
\begin{enumerate}
\item $\left.\begin{array}{rcl}
            x(t+\om) &=& x(t) + x^* \\
             a(t+\om) &=& a(t) + a^* + \frac{1}{2}\left[x^*,x(t)\right]
             \end{array}\right\}$ for all $t\in\R$ and some $\om>0$;

\item $\g(t+\om) = \phi\g(t)$ for all $t\in\R$ and some $\om>0$;

\item $e^{\om J}$ fixes $e_1 + y_1 + x_2 = e_0 + y_1 + J^{-1}y_2$.
\end{enumerate}
\end{proposition}
\begin{proof}
As before, we may assume without loss of generality  $n = 1\in N$.
Items 1 and 2 are equivalent by formula (\ref{bch}).  The following lemma
shows that item 1 implies item 3.
\begin{lemma}
As\label{lem} in the preamble to Theorem \ref{ige}, assume $\E$ as an
orthogonal direct sum $\E_1\ds\E_2$ with $\E_1 = \ker J$, and use $x_1$,
$x_2$, $y_1$, and $y_2$ as given there. Then $x^* = \om x_1 + \half\om^2
y_1$ and $e^{\om J}$ fixes $x_2$.
\end{lemma}
\begin{proof}
By Theorem \ref{ige}, $x(k\om) = k\om x_1 + \left(e^{k\om J} -
I\right) J^{-1}x_2 + \half k^2\om^2 y_1$ for every positive integer $k$.
By induction, from item 1 in the statement of the proposition we obtain
$x(k\om) = k x^*$ for every $k$. One may decompose $x^* = v^* + e^*_1 + e^*_2
\in \V\ds\E_1\ds\E_2$ to see $v^* = \om v_0$, $e^*_1 = \om e_1 +
\half\om^2 y_1$, and
\begin{equation}
k( e^*_2 + \om J^{-1}y_2) = \left(e^{k\om J} - I\right)J^{-1}x_2\ \forall\ k\,.
\label{unbdd}
\end{equation}
for every $k$.

Now, $e^{\om J}$ is an element of the identity component of the
pseudorthogonal group of isometries of $\langle\,,\rangle$, and as such
can be decomposed into a product of reflections, ordinary rotations, and
boosts. With respect to appropriate coordinates, which may be different
from our standard choice, a boost will have a matrix of the form
$$ \left[\begin{array}{cc}
         \cosh s & \sinh s\\
         \sinh s & \cosh s \end{array}\right] $$
on some pair of basis vectors, for some $s\in\R$.

If $e^{\om J}$ is composed only of reflections and ordinary rotations,
then the right-hand side of (\ref{unbdd}) is uniformly bounded in $k$
(say, with respect to the positive definite $\langle\,,\io\rangle$) while
the left-hand side is unbounded, we see that $e^*_2 + \om J^{-1}y_2 = 0$.  On the
other hand, if $e^{\om J}$ is a pure boost, then the right-hand side grows
exponentially in $k$ while the left-hand side grows but linearly, and
again we obtain $e^*_2 + \om J^{-1}y_2 = 0$.  The Lemma now follows from
this and (\ref{unbdd}) for $k=1$.
\end{proof}
The proof that item 3 implies item 2 is the same as the relevant part of
the proof of (4.3) in \cite{E1}.
\end{proof}
\begin{corollary}
When\label{gctgc} in addition $v^* = v_0 = 0$, the following are also
equivalent to the three items in Proposition \ref{gctg}.
\begin{enumerate}
\item $\dg(0)$ is orthogonal to the orbit $Z_{e^*}n$, where
$Z_{e^*} = \exp\left([e^*,\n]\right) \subseteq Z$;

\item $\dg(\om)$ is orthogonal to the orbit $Z_{e^*}\phi n$.
\end{enumerate}
\end{corollary}
\begin{proof}
Note that under this hypothesis, $x^* = e^*$.  Now Lemma \ref{lem}
implies that $J(e^*) = 0$, and this is now equivalent to $z_0 \perp
[e^*,\n]$.  Thus the relevant parts of the proof of (4.3) in \cite{E1}
apply {\em mutatis mutandis.}
\end{proof}

We also obtain the following results as in Eberlein
\cite[(4.4),\,(4.9)]{E1}.  Note that we assume that $v^* = v_0 = 0$ in the
first, but that this is automatic in the second.
\begin{corollary}
Let\label{tsg1} $\phi\in N$ and write $\phi = \exp(a^* + e^*)$ for unique
elements $a^*\in\z$ and $e^*\in\E$. Let $n\in N$ be given and write $n =
\exp(\xi)$ for a unique $\xi\in\n$. Then the following are equivalent:
\begin{enumerate}
\item There exists a geodesic $\g$ in $N$ with $\g(0) = n$ such
that $\phi\g(t) = \g(t + \om)$ for all $t\in\R$ and some $\om > 0$.

\item There exists a geodesic $\g^*$ in $N$ with $\g^*(0) = 1$,
$\dg^*(0)$ is orthogonal to $[e^*,\n]$, and $\g^*(\om)
= \exp\left([e^*, \xi]\right)\phi$ for some $\om > 0$.\eop
\end{enumerate}
\end{corollary}
\begin{corollary}
Let\label{tsg2} $1\ne\phi\in Z$ and $\g$ be any geodesic so $\g(\om)
= \phi\g(0)$ for some $\om > 0$.  Then $\phi\g(t) = \g(t + \om)$ for all
$t\in\R$.\eop
\end{corollary}

In the flat 2-step nilmanifolds of Theorem \ref{e0f}, we can calculate the
period spectrum completely.
\begin{theorem}
If\label{upd} $[\n,\n] \subseteq \U$ and $\E=\{0\}$, then $\specp(M)$ can be
completely calculated from $\log\G$ for any $M = \G\bs N$.
\end{theorem}
\begin{proof}
Let $\phi$ translate a unit-speed geodesic $\g$ by $\om > 0$.  As usual,
we may as well assume that $\g(0) = 1 \in N$.  Write $\log\phi = a^* +
v^*$ and $\dg(0) = a_0 + v_0$.  From Corollary \ref{iged} we have that $v^* =
v(\om) = \om\,v_0$, $z^* = z(\om) = \om\,z_0$, and $u^* = u(\om) =
\om\,u_0 + \half\om^2\,\sJ v_0$.  Note that $$\om^2\,\sJ v_0 =
\om^2\,\add{v_0}{(z_0 + v_0)} = \om^2\,\add{v_0}{v_0} = \add{v^*}{v^*}\,.$$
Substituting and rearranging, we obtain
$ \ve\om^2 = 2\langle u^*,v^*\rangle + \langle z^*,z^*\rangle$,
where $\pm1 = \ve = \langle\dg(0),\dg(0)\rangle = 2\langle u_0,v_0\rangle
+ \langle z_0,z_0\rangle$.
\end{proof}
Thus we see again, as mentioned after Corollary \ref{iged}, just how much
these flat, 2-step nilmanifolds are like tori.  All periods can be
calculated purely from $\log\G \subseteq \n$, although some will not show
up from the tori in the fibration.
\begin{corollary}
$\specp(T_B)$ (respectively, $T_F$) is\label{psd} $\cup_{\sC}\,\wp(\sC)$
where the union is taken over all those free homotopy classes\/ $\sC$ of
closed curves in $M = \G\bs N$ that\/ {\em do not} (respectively,\/ {\em
do}) contain an element in the center of\/ $\G \cong \pi_1(M)$, except for
those periods arising only from unit-speed geodesics in $M$ that project to
null geodesics in both $T_B$ and $T_F$.\eop
\end{corollary}
We note that one might consider using this to assign periods to some null
geodesics in the tori $T_B$ and $T_F$.

When the center is nondegenerate, we obtain results similar to Eberlein's
\cite[(4.5)]{E1}.
\begin{proposition}
Assume\label{e4.5} $\U = \{0\}$. Let $\phi \in N$ and write $\log\phi =
z^* + e^*$.  Assume $\phi$ translates the unit-speed geodesic $\g$ by $\om
> 0$.  Let $z'$ denote the component of $z^*$ orthogonal to $[e^*,\n]$.
Let $n = \g(0)$ and set $\om^* = |z' + e^*|$. Let $\dg(0) = L_{n*}(z_0 +
e_0)$ and use $J$, $e_1$, and $e_2$ as in Corollary \ref{igend} (see also
just before Theorem \ref{ige}).  Then
\begin{enumerate}
\item $|e^*| \le \om$. In addition, $\om < \om^*$ for timelike
   (spacelike) geodesics with $\om z_0 - z'$ timelike (spacelike), and
   $\om > \om^*$ for timelike (spacelike) geodesics with $\om z_0 - z'$
   spacelike (timelike).

\item $\om = |e^*|$ if and only if $\g(t) = \exp\left(t\,e^*\!/|e^*|
   \right)$ for all $t \in \R$.

\item $\om = \om^*$ if and only if $\om z_0 - z'$ is null.  If moreover
   $\om^* z_0 = z'$, then $e_2 = 0$ if and only if
\begin{enumerate}
   \item $\g(t) = n\,\exp\left( t\,(z' + e^*)/\om^*\right)$ for all
      $t \in \R$.

   \item $z' = z^* + [e^*,\xi]$ where $\xi = \log n$.
\end{enumerate}
\end{enumerate}
\end{proposition}
Although $\om^*$ need not be an upper bound for periods as in the
Riemannian case, it nonetheless plays a special role among all periods,
as seen in item 3 above, and we shall refer to it as the {\em
distinguished\/} period associated with $\phi \in N$.  When the center is
definite, for example, we do have $\om \le \om^*$.
\begin{proof}
As usual, we may assume that $\g(0) = 1 \in N$.  Note that this replaces
$\phi$ as given in the statement with $n^{-1}\phi n$ and $\g$ with
$n^{-1}\g$.

For the first part of item 1, since $|\dg(0)| = 1$ there exists an
orthonormal basis of $\n$ having $\dg(0)$ as a member.  (This may well be
a different basis from our usual one.)  Fix one such basis, and consider
the positive-definite inner product with matrix $I$ on this basis.  Let
$\|\cdot\|$ denote the norm associated to this positive-definite inner
product. By Lemma \ref{lem}, $e^* = \om e_1$.  Then $|e_1| \le \|e_1\| \le
\|\dg(0)\| = 1$ so $|e^*| =
\om |e_1| \le \om$.

For the rest of item 1, we begin with Corollary \ref{igend} and get
\begin{eqnarray*}
e(t) &=& t\,e_1 + \left(e^{tJ} - I\right)J^{-2}e_2\,,\\
z(t) &=& t\left( z_0 + \half\left[ e_1, \left(e^{tJ}+I\right)J^{-1}e_2
         \right] \right) + z_2(t) + z_3(t)\,.
\end{eqnarray*}
By Lemma \ref{lem}, $e^* = \om e_1$ and $e^{\om J}e_2 = e_2$. Inspecting
the formula for $z_2(t)$ in Corollary \ref{igend}, we find $z_2(\om) = 0$.
Thus
\begin{eqnarray*}
z^* = z(\om) &=& \om\left(z_0 + [e_1,J^{-1}e_2]\right) + z_3(\om)\\
   &=& \om\,z_0 + [e^*,J^{-1}e_2] + \half\int_0^{\om} \left[ e^{sJ}
J^{-1}e_2, e^{sJ}e_2\right] ds\,.
\end{eqnarray*}
By item 1 of Corollary \ref{gctgc}, $z_0 \perp [e^*,\n]$.  Then
$$ \langle z',z_0\rangle = \langle z^*,z_0\rangle = \om\langle
z_0,z_0\rangle + \half\int_0^{\om}\left\langle\left[e^{sJ}J^{-1}e_2,
e^{sJ}e_2\right], z_0\right\rangle ds\,.$$
Recall that $J$ is skewadjoint with respect to $\langle\,,\rangle$ (whence
$e^{tJ}$ is an isometry of $\langle\,,\rangle$ for all $t$), that $J$
commutes with every $e^{tJ}$ (whence so does $J^{-1}$), and that $Jx =
\add{x}{z_0}$. We compute
\begin{eqnarray*}
\left\langle\left[e^{sJ}J^{-1}e_2, e^{sJ}e_2\right], z_0\right\rangle
&=& -\left\langle\left[e^{sJ}e_2, e^{sJ}J^{-1}e_2\right], z_0\right\rangle\\
&=& -\left\langle J^{-1}e^{sJ}e_2, Je^{sJ}e_2\right\rangle\\
&=& \left\langle e^{sJ}e_2, e^{sJ}e_2\right\rangle\\
&=& \langle e_2, e_2\rangle.
\end{eqnarray*}
Therefore,
\begin{equation}
\langle z', z_0\rangle = \om\langle z_0, z_0\rangle +
\frac{\om}{2}\langle e_2, e_2\rangle .\label{one}
\end{equation}

Now $|\dg(0)| = 1$ so $\ve = \langle z_0,z_0\rangle + \langle e_1,
e_1\rangle + \langle e_2,e_2\rangle$, where $\ve = \pm 1$ as usual.
Substituting in (\ref{one}) for $\langle e_2,e_2\rangle$, we obtain
$$ \langle z',z_0\rangle = \frac{\om}{2}\bigl( \ve + \langle
   z_0,z_0\rangle\bigr) - \frac{\om}{2}\frac{\langle
   e^*,e^*\rangle}{\om^2} $$
so
$\langle e^*,e^*\rangle - \ve\om^2 = \om^2\langle z_0,z_0\rangle
-2\om\langle z',z_0\rangle$.
Adding $\langle z',z'\rangle$ to both sides,
\begin{equation}
\langle z'+e^*, z'+e^*\rangle - \ve\om^2 =
   \langle \om z_0-z', \om z_0-z'\rangle.\label{two}
\end{equation}
There are several cases: $\ve$ is 1 or $-1$ and $\om z_0-z'$ is timelike,
spacelike, or null.  If $\om z_0-z'$ is null, then $|\langle z'+e^*,
z'+e^*\rangle| = \om^2 > 0$ and $\om = |z'+e^*|$.  If $\ve = 1$ and $\om
z_0-z'$ is timelike, or if $\ve = -1$ and $\om z_0 - z'$ is spacelike,
then $\ve\langle z'+e^*, z'+e^*\rangle > \om^2 > 0$ whence $\om <
|z'+e^*|$.  If $\ve = 1$ and $\om z_0-z'$ is spacelike, or if $\ve = -1$
and $\om z_0 - z'$ is timelike, then it follows similarly that $\om > |z'
+ e^*|$.  This completes the proof of item 1.

Now we prove item 2.  If $\g$ is as given there, then $\om = |e^*|$ because $\exp(z^* + e^*) =
\phi = \g(\om) = \exp(\om\,e^*\!/|e^*|)$.
Conversely, assume $\om = |e^*|$ and consider the associated
positive-definite inner product $\langle\cdot,\io\cdot\rangle$.  Changing
the basis of $\E$ if necessary, we may assume that $\Z$, $\E_1$, and
$\E_2$ are mutually orthogonal with respect to both $\langle\,,\rangle$
and $\langle\cdot,\io\cdot\rangle$.  Let $\|\cdot\|$ now denote the norm
for $\langle\cdot,\io\cdot\rangle$.  Then
\begin{equation}
\|\dg(0)\|^2 = \|z_0\|^2 + \|e_1\|^2 + \|e_2\|^2\label{three}
\end{equation}
so $\|\dg(0)\|^2 = \|e_1\|^2$ if and only if $\dg(0) = e_1 = e^*\!/|e^*|$.
But now $\g$ has the same initial data as $\exp(t\,e^*\!/|e^*|)$, so by
uniqueness they must coincide.

Finally, we prove the last part of item 3; the first part is immediate
from the last part of the proof of item 1 above.  So assume $\om^*z_0 - z'
= 0$ or $z_0 = z'\!/\om^*$. Continue with the immediately previous
positive-definite norm $\|\cdot\|$ and basis of $\E$.  Substituting in
(\ref{three}) we get
\begin{eqnarray*}
\|\dg(0)\|^2 &=& \frac{\|z'\|^2}{(\om^*)^2} + \frac{\|e^*\|^2}{(\om^*)^2}
+ \|e_2\|^2=\frac{\|z'+e^*\|^2}{(\om^*)^2} + \|e_2\|^2
\end{eqnarray*}
whence
$\dg(0) = \frac{z'+e^*}{\om^*}$
if and only if $e_2 = 0$.
\end{proof}
\begin{corollary}
Assume\label{e4.6} the center is nondegenerate.  Let $\phi \in N$ with
$\phi \notin Z$ and suppose that $z^* \in [e^*,\n]$. Then
\begin{enumerate}
\item If $\phi$ translates a timelike (spacelike) geodesic with $z_0$
nonspacelike (nontimelike), then $\phi$ has the unique period $|e^*|$.

\item Let $\g$ be a unit-speed geodesic in $N$ with $\g(0) = n =
\exp(\xi)$ for a unique $\xi \in \n$. Then $\phi$ translates $\g$ by the
unique period $|e^*| > 0$ if and only if $[\xi,e^*] = z^*$ and $\g(t) =
n\,\exp(t\,e^*\!/|e^*|)$ for all $t \in \R$.
\end{enumerate}
In particular, this applies to all noncentral $\phi \in N$ if $\n$ is
nonsingular.\eop
\end{corollary}
The proof follows that of \cite[(4.6)]{E1} {\em mutatis mutandis\/} and we
omit the details. From item 2 of Proposition \ref{e4.5},
using Lemma \ref{dexp}, we obtain 
\begin{eqnarray*}&&\exp(e^* + \half[\xi,e^*]) =
n\,\exp(e^*) = \g(|e^*|) = \phi\,n \\&=& \exp(z^* + e^*)\exp(\xi) = \exp(z^* +
e^* + \half[e^*,\xi]),\end{eqnarray*} thus avoiding the use of item 3 here.  Anent the
last comment, note that if $\n$ is nonsingular then in fact $z_0 = 0$ in
item 1, because $z_0 \perp [e^*,\n] = \z$.

In view of the comment following Proposition \ref{e4.5}
and Corollary \ref{e4.6}, the following definitions make sense at
least for $N$ with a nondegenerate center.
\begin{definition}
Let\label{MpC} $\sC$ denote either a nontrivial, free homotopy class of
closed curves in $\G\bs N$ or the corresponding conjugacy class in
$\G$. We define $\wp^*(\sC)$ to be the distinguished periods of periodic
unit-speed geodesics that belong to $\sC$.
\end{definition}
\begin{definition}
The\label{Dpsp} {\em distinguished period spectrum\/} of $\G\bs N$ is 
the set
$$ \Dspecp(\G\bs N) = \bigcup_{\sC}\wp^*(\sC)\,,$$
where the union is taken over all nontrivial, free homotopy classes of
closed curves in $\G\bs N$.
\end{definition}
Then as an immediate consequence of the preceding corollary, we get:
\begin{corollary}
Assume\label{psnd} the center is nondegenerate.  If $\n$ is nonsingular,
then $\specp(T_B)$ (respectively, $T_F$) is precisely the period spectrum
(respectively, the distinguished period spectrum) of those free homotopy
classes\/ $\sC$ of closed curves in $M = \G\bs N$ that\/ {\em do not}
(respectively,\/ {\em do}) contain an element in the center of\/ $\G \cong
\pi_1(M)$, except for those periods arising only from unit-speed geodesics
in $M$ that project to null geodesics in both $T_B$ and $T_F$.\eop
\end{corollary}

\section*{Acknowledgments}
Once again, Parker wishes to thank the Departamento at Santiago for its
fine hospitality.  He also thanks WSU for a Summer Research Fellowship
during which part of this work was done, and for a Sabbatical Leave during
which it was continued.

\end{document}